\RequirePackage{times}
\RequirePackage{fullpage}
\RequirePackage{ifthen}

\documentclass{article} 

\usepackage{titlesec}
\usepackage[a4paper, total={7.5in, 10in}]{geometry}
\usepackage{bm}
\usepackage[latin1]{inputenc}
\usepackage{amsmath}
\usepackage{amsfonts}
\usepackage{amssymb}
\usepackage{graphicx}
\usepackage{notoccite}
\usepackage{longtable}
\usepackage[numbers,sort&compress]{natbib}

\titleformat*{\section}{\bfseries}
\titleformat*{\subsection}{\bfseries}

\newcommand{\spacing}[1]{\renewcommand{\baselinestretch}{#1}\large\normalsize}
\spacing{1}

\begin{document}

\newcommand{\Mu}{\mathcal{M}}
\newcommand{\I}{\mathcal{I}}

	\begin{center}
		{\Large		
			\textbf\newline{\bf \textsf{Bayesian information 
					criteria for clustering 
					normally distributed data}}
		}
		\vspace{0.2cm}
		\newline
		Anthony J. Webster
		\\
		\vspace{0.15cm}
		{\sl 
		Nuffield Department of Population Health, 
		Big Data Institute, Old Road Campus, University of Oxford, 
		Oxford, OX3 7LF,
		UK.}
		\\
	\end{center}

\begin{abstract}
{\bf 
Maximum likelihood estimates (MLEs) are asymptotically normally distributed, 
and this property is used in meta-analyses to test the heterogeneity of 
estimates, either for a single cluster or for several sub-groups.  
More recently, MLEs for associations between risk factors and diseases have been 
hierarchically clustered to search for diseases with shared underlying causes, 
but an objective statistical criterion is needed to determine the  
number and composition of clusters.  
To tackle this problem, 
conventional statistical tests are briefly reviewed, before considering the 
posterior distribution for a partition of data into clusters. 
The posterior distribution is calculated by marginalising out the unknown 
cluster centres, and is different to the likelihood associated with 
mixture models.  
The calculation is equivalent to that used to obtain the Bayesian Information 
Criterion (BIC), but is exact, without a Laplace approximation.
The result includes a sum of squares term, and terms that depend on the 
number and composition of clusters, that penalise the number of free parameters 
in the model.
The usual BIC is shown to be unsuitable for clustering applications unless the number of items in each individual cluster is sufficiently large. 
}
\end{abstract}

\section{Introduction}

Despite some interest in clustering data sampled from normal distributions 
\cite{Webster2021,Nielsen2009,Teklehaymanot2018}, most recent work has 
focussed on clustering distributions of data  that are non-normal or containing outliers  \cite{Punzo2020,Dotto2019,Bagnato2017,Lin2014,Lee2013,Contreras2012}, 
or has focussed on developing methods for clustering data whose underlying 
distributions are unknown  \cite{Wasserman2010,Bishop2006}. 
However, multivariate normally distributed data commonly arise, 
in particular as the asymptotic distribution of maximum 
likelihood estimates (MLEs). 
When clustering MLEs, their (estimated) covariances are given and must be accounted for, and these are determined by 
sample sizes and often strong correlations between covariates. 
Here the posterior probability for clusters of normally distributed data with 
known covariances is calculated by marginalising out the unknown cluster centres. 
The result is equivalent to exactly calculating the Bayesian Information Criterion 
(BIC) \cite{Wasserman2010,Schwarz1978,Djuric1998,Cavanaugh1999,Teklehaymanot2018}, 
without the Laplace approximation.  
It includes a weighted sum of squares term to penalise 
poor fits, and a term equivalent to the $k\log n$ term in BIC  \cite{Wasserman2010,Schwarz1978,Djuric1998,Cavanaugh1999,Teklehaymanot2018}, 
where $k$ and $n$ are the number of parameters and data respectively.  
For recent work on approximate BIC calculations for use with clustering, 
and a review of related literature, see \cite{Teklehaymanot2018}. 
The greatest difference to other studies, is that the work 
here assumes that each data point is sampled from a normal distribution 
with a known covariance and a mean that is an (unknown) cluster 
centre. 
When the underlying distribution of data are unknown, there are advantages to 
firstly fitting a parametric model, and then clustering the estimated co-efficients, 
but the approach is rarely used  \cite{Holmes2005,Heard2006,Kim2008,Kim2017,Park2019,Teklehaymanot2018}. 
A recent example \cite{Webster2021} used this approach to search for shared 
underlying causes of disease, in a large study using UK Biobank data \cite{UKBiobank}. 
In that case, existing epidemiological understanding was incorporated through 
the selection, truncation, and censoring of data, and through the adjustment 
for known confounding factors with established survival analysis methods \cite{Webster2021}. 
Advantages of the general approach include: 
\begin{enumerate}
	
	\item{The distribution of underlying data can be unknown, but 
		under reasonable regularity conditions, MLE estimates 
		are normally distributed \cite{Wasserman2010,Hardle2015}.}
	
	\item{Estimates can be stratified and multiply-adjusted 
		for known confounders, helping to extract the information of interest 
		from potentially noisy 
		data \cite{Webster2021,Wasserman2010,Collett2014}.}
	
	\item{The fitted models can be more interpretable and familiar 
		to the scientific community. For example, 
		proportional hazards models are commonly used by medical researchers.}
	
	\item{Marginalisation \cite{Wasserman2010,Bishop2006} can be used 
		to select parameter subsets of most interest for 
		clustering \cite{Webster2021}, 
		e.g. risk factors as 
		opposed to confounders, or minimum versus maximum quantiles.}
	
	\item{By fitting a model, we can build-in
		prior knowledge through the model.}
	
\end{enumerate}
In addition to clustering of diseases \cite{Webster2021}, the approach has 
been used to detect changes in 
gene expression by clustering Fourier series coefficients \cite{Kim2008,Kim2017}. 
Clustering parameters from linear-models such as a Fourier series, are examples of 
the more general problem of clustering normally-distributed MLE estimates  
from a parametric model. 
Here we consider the general problem of clustering data sampled from multivariate 
normal distributions, with the aim of determining the number and membership of clusters. 
A related approach for clustering categorical data with an 
``exact integrated complete-data likelihood (ICL)'', was developed 
by Biernacki et al. \cite{Biernacki2000,Biernacki2010}, and 
has been used for model selection in several 
applications \cite{Rigaill2012, Come2015, Marbac2017, Lomet2018, Come2021}, 
including a recent clustering of human population genomic 
data \cite{Marbac2020}. 
An important difference between the model described here, the ICL, and 
clustering with Gaussian mixture models for example, is that the 
MLE data studied here each have an (asymptotically) normal distribution 
with an estimated mean and covariance. 
As a result, the probability distributions are quite 
different, despite all three approaches having a wide range of 
applicability. 
The next section briefly reviews statistical tests of heterogeneity, 
their use in meta-analyses, and their potential use in clustering studies. 
The formulae involve sums of squares, and similar terms occur in 
Section \ref{posterior}, that calculates the posterior distribution associated 
with a partition of data into clusters. 
Section \ref{bic} relates the results to the Bayesian information 
criterion (BIC), and indicates the conditions when the usual BIC could be used in 
clustering studies, and how this would be done. 
It also emphasises that in general, the usual BIC cannot be used in 
clustering studies. 
The final Section \ref{discussion} discusses: limits of the prior distribution, 
similarities to clustering with k-means, potential uses in meta-analyses, the 
scope for sensitivity analyses or forming confidence sets, and potential model 
improvements. 

\section{Statistical heterogeneity tests}\label{statTests}

Heterogeneity tests are widely used in meta-analyses, and are intended to 
assess whether estimates are the same in several 
different studies \cite{Borenstein2009}. 
These include multivariate heterogeneity tests  \cite{Kalaian1996,vH2002,Nam2003,Jackson2012}, 
and the tests can involve both fixed effects  
and random effects models \cite{Borenstein2009}. 
A fixed effects model is considered first, 
that generalises easily to a random effects model. 
In a fixed effects model the null hypothesis is that all diseases 
have the same associations with one or more 
parameters, such as a drug, or a collection 
of potential risk factors.
These might be a subset of associations, with potential confounders 
adjusted for, and subsequently removed by marginalisation \cite{Webster2021}. 
Consider $m$ items, such as a collection of diseases, in a cluster labelled by $g$. 
Under the null hypothesis the $i$th item will have, 
\begin{equation}\label{Xi}
	X_i \sim N\left( \mu_g, \Gamma_i^{-1} \right) 
\end{equation}
where $\Gamma_i^{-1}=\Sigma_i$ is the covariance 
($\Gamma_i$ is the precision matrix), and $\mu_g$ are an (unknown) vector of 
the estimated associations, that are assumed to be the same for 
all items in the cluster (e.g. diseases in a composite endpoint). 
Eq. \ref{Xi} requires \cite{Hardle}, 
\begin{equation}
	\left(X_i-\mu_g\right)^T \Gamma_i \left( X_i-\mu_g \right) \sim \chi^2_p
\end{equation}   
where $p$ is the dimension. 
Therefore because the sum of $n_g$ random variables that are individually 
$\chi^2_p$ distributed is $\chi^2_{pn_g}$, 
\begin{equation}\label{chisqXi}
	\sum_{i\in C_g}  
	\left(X_i-\mu_g\right)^T \Gamma_i \left( X_i-\mu_g \right) \sim \chi^2_{p n_g}
\end{equation}
where $C_g$ is the set of $n_g$ items in cluster $g$. 
For $p=1$ this has, 
\begin{equation}
	\sum_{i\in C_g} \left(\frac{X_i-\mu_g}{\sigma_i} \right)^2 \sim \chi^2_{n_g}
\end{equation}
for standard deviations $\sigma_i$. 
Because $\mu_g$ is unknown and must be estimated, the test statistic is 
modified, as explained next. 
Appendix \ref{BayesApp} uses Bayes theorem with a flat or normal prior for the 
mean $\mu_g$, to show that if $X_i \sim N\left( \mu_g, \Gamma_i^{-1} \right)$, then, 
\begin{equation}\label{Pmug}
	\mu_g \sim N\left( \tilde{\mu}_g, \tilde{\Gamma}_g^{-1} \right) 
\end{equation}
where,
\begin{equation}\label{BmuG}
	\tilde{\mu}_g =   
	\left( \sum_{i} \Gamma_i \right)^{-1} 
	\left( \sum_{i} \Gamma_i X_i \right)
\end{equation} 
and, 
\begin{equation}\label{LambdaG}
	\tilde{\Gamma}_g = \left( \sum_{i} \Gamma_i \right)
\end{equation}
If a normal prior is used then the sum over $i$ includes the prior's 
mean $\mu_0$ and covariance $\Gamma_0$, and the sum is from $i=0$ to $i=n_g$. 
For a flat prior, the sum is from $i=1$ to $i=n_g$.
The subscripts $g$ will later allow the discussion to include more than one cluster, 
for example several clusters of diseases as were considered in 
Webster et al. \cite{Webster2021}. 
For the rest of this section, unless stated otherwise, we consider a single cluster. 
For $p=1$, Eq. \ref{BmuG} is the well-known inverse-variance weighted 
estimate of the mean. 
Appendix A shows that,
\begin{equation}\label{penult}
	\begin{array}{c}
		\sum_i 
		\left( X_i - \tilde{\mu}_g  \right)^T 
		\Gamma_i 
		\left( X_i - \tilde{\mu}_g  \right)
		\\
		=
		\sum_i \left( X_i - \mu_g \right)^T \Gamma_i \left( X_i -\mu_g \right)
		-  
		\left(  \mu_g - \tilde{\mu}_g  \right)^T 
		\tilde{\Gamma}_g 
		\left(  \mu_g - \tilde{\mu}_g  \right)
	\end{array}
\end{equation}
and Eq. \ref{Pmug}, implies that, 
\begin{equation}\label{chisqmu}
	\left( \mu_g -\tilde{\mu}_g \right)^T \tilde{\Gamma}_g 
	\left( \mu_g -\tilde{\mu}_g \right) \sim \chi^2_p 
\end{equation}
These observations can be used to 
test the assumption that the normal distributions have the same mean. 
Using  Eqs. \ref{chisqXi}. \ref{penult}, and \ref{chisqmu}, 
\begin{equation}\label{cepss}
	\sum_{i\in C_g}
	\left(  X_i - \tilde{\mu}_g  \right)^T 
	\Gamma_i 
	\left(  X_i - \tilde{\mu}_g \right)
	\sim \chi^2_{(n_g-1)p}
\end{equation}
The left side of Eq. \ref{cepss} is the Q statistic. 
It tests the assumption that a set of  
(approximately) normally distributed estimates $\{X_i\}$ have the 
same mean. 
For $p=1$, these expressions give the well known inverse 
variance weighted heterogeneity test, that is regularly 
used in meta analyses and 2-sample Mendelian 
randomisation studies \cite{Borenstein2009,Burgess2014}.
For the situation described in Webster et al. \cite{Webster2021}, the aim is to 
assess the goodness of fit for a clustering of diseases. 
For this situation, Eq. \ref{cepss} is modified to sum over all $m$ clusters, 
and the Q statistic becomes, 
\begin{equation}\label{clusteringQsq}
	\sum_{g=1}^m
	\sum_{i\in C_g} 
	\left(  X_i - \tilde{\mu}_g  \right)^T 
	\Gamma_i 
	\left(  X_i - \tilde{\mu}_g \right)
	\sim   
	\chi^2_{(n-m)p}
\end{equation} 
where we used $\sum_{g=1}^m(n_g-1)p=(n-m)p$, and 
$C_g$ is the set of diseases in cluster $g$ (composite endpoint $g$). 
If there is a normal prior for the cluster centres, as discussed in 
Sections \ref{posterior} and \ref{bic}, then $n_g$ will 
include one extra element for each $g$.
As a consequence, the right sides of Eqs \ref{cepss} and \ref{clusteringQsq} 
become, $\sim \chi^2_{n_g p}$ and $\sim \chi^2_{n p}$ respectively. 
The results above correspond to a fixed effects model where the data are 
assumed to have the same mean, as opposed to the means being sampled from 
an underlying distribution (a random effects model).
A random effects model firstly assumes that each study has a mean $\mu_g$ with  
$\mu_g \sim N(\mu_{g_0}, \Gamma_{g_0}^{-1})$, but that the measured 
estimates $X_i$ have 
$X_i\sim N(\mu_g,\Gamma_i^{-1})$. 
It then marginalises over $\mu_g$ for each $i$, to give  
$X_i \sim N(\mu_{g_0}, \Gamma_{g_0}^{-1}+\Gamma_i^{-1})$, 
that replaces Eq. \ref{Xi}.
Eqs. \ref{BmuG} and 
\ref{LambdaG} become, 
\begin{equation}
	\tilde{\mu}_g = \left(
	\sum_i \left( \Gamma_{g_0}^{-1} + \Gamma_i^{-1} \right)^{-1} 
	\right)^{-1}
	\left(
	\sum_i \left( \Gamma_{g_0}^{-1} + \Gamma_i^{-1} \right)^{-1} X_i 
	\right)
\end{equation}
and
\begin{equation}
	\tilde{\Gamma}_{g_0} = 
	\sum_i \left( \Gamma_{g_0}^{-1} + \Gamma_i^{-1} \right)^{-1} 
\end{equation}
as in Jackson et al. \cite{Jackson2012} (where $\tilde{\Lambda}_0$ 
is the inverse of the covariance). 
To calculate the $Q$ and $I^2$ statistics, 
$\Gamma_i$ is replaced by $(\Gamma_{g_0}^{-1}+\Gamma_i^{-1})^{-1}$, 
and $n_g$ is the total number of studies.   
The above arguments and results could be modified to consider a random effects 
model with different priors for each cluster. 
The random effects model is equivalent to treating each data 
point as an individual cluster whose centre is sampled from a normal distribution. 
This contrasts with the cluster model with a normal prior, where the clusters' 
centres are only sampled once per cluster. 
The latter is equivalent to a subgroup analysis in which (normally distributed) 
random heterogeneity is assumed between the subgroups, but 
a fixed effects model is assumed for the heterogeneity within groups. 

\subsection{Heterogeneity and meta-analyses}\label{hetTests}

The traditional measure of heterogeneity is the Q statistic, that was 
derived above in a multi-variate context.
The I-square statistic \cite{Higgins2002,Higgins2003} is closely related 
to the Q statistic \cite{Borenstein2009}, and in the 
notation above, is for Eq. \ref{clusteringQsq},
\begin{equation}\label{Isq}
I^2 = \frac{Q - (n-m)p}{Q} \times 100
\end{equation}
where $Q$ is the left side of Eq. \ref{clusteringQsq}, and the 
factor of $100$ is conventionally used to express $I^2$ as a percentage. 
$I^2$ is usually set to zero if its evaluation is negative. 
The equivalent expression for a single cluster that uses the left side 
of Eq. \ref{cepss} for $Q$, would replace the number of degrees of 
freedom $(n-m)p$, with $(n_g-1)p$ in Eq. \ref{Isq}. 
If as discussed in Sections \ref{posterior} and \ref{bic}, there is a normal 
prior for the cluster centres,  then $n_g$ will 
include one extra element for each $g$, causing $(n-m)p$ to be replaced 
by $np$, and $(n_g-1)p$ to be replaced by $n_g p$. 
The $I^2$ statistic replaces a test with a more nuanced measure of 
heterogeneity that is useful when some heterogeneity is 
expected, but it does not provide an objective statistical test. 

\section{Posterior distribution for clusters of normally distributed data}\label{posterior}

Firstly consider a distribution \cite{Marcus} of identifiable 
items in identifiable boxes. Later this will be used this to consider 
the likelihood for a partition \cite{Marcus} of identifiable  
items in unlabelled clusters\footnote{Implicitly, a partition's clusters 
	are labelled by the number of items they contain, but this does not uniquely 
	label the clusters in a partition because more than one cluster can have the same number of elements.}. 
{Notation:} The $i$th item's cluster 
in a {distribution} of labelled clusters, is $Z_i$. 
The MLEs $\{ \hat{\mu}_i\}$  and their covariances 
$\{ \hat{\Sigma}_i \}$, asymptotically have, 
\begin{equation}\label{muEq}
	\hat{\mu}_i \sim N \left( \mu_{Z_i}, \hat{\Sigma}_i \right) 
\end{equation}
with $\mu_{Z_i}=\mu_{Z_j}$ iff $Z_i=Z_j$.
We will regard $\hat{\Sigma}_i$ as given, and $\hat{\mu}_i$ as 
random variables sampled from Eq. \ref{muEq}.  
Write $X_i = \hat{\mu}_i$,  
$\Gamma_i = \hat{\Sigma}_i^{-1}$,
and the propositions 
$Z=\{Z_i=z_i\}$,
$X=\{X_i=x_i\}$, 
$\Mu=\{\Mu_g=\mu_g\}$, where $\mu_g$ is the mean of cluster $g$.
We will also write 
$\I$ to denote both 
$\{\Gamma_i\}$ and any additional implicit information regarding 
the clustering problem. 
Note that the proposition $Z=\{Z_i=z_i\}$ is not simply a 
conjunction of independent terms.
For example, $Z$ implicitly contains 
information about the number of clusters $M=m$, the number of members  
in each cluster, and the number of $k$-element clusters. 
One consequence is that we cannot 
directly write $P(Z|X,\I)$ as a product 
with one term per item to be clustered. 
Consider the likelihood of assigning labelled items into 
a distribution \cite{Marcus} of labelled clusters. 
Using Bayes theorem and marginalisation \cite{Wasserman2010,Bishop2006},  
\begin{equation}\label{bayes1}
	\begin{array}{ll}
		P \left( Z | X, \I \right)P(X|\I)
		&=
		P \left( X | Z, \I \right) 
		P\left( Z  | \I \right)
		\\
		&=
		\int_{-\infty}^{\infty} d\mu_1 ... \int_{-\infty}^{\infty} d\mu_m
		P \left( X, \Mu |  Z, \I \right) 
		P\left( Z | \I \right) 
		\\
		&=
		\int_{-\infty}^{\infty} d\mu_1 ... \int_{-\infty}^{\infty} d\mu_m
		P \left( X | \Mu,  Z, \I \right) 
		P \left( \Mu | Z, \I \right) 
		P\left( Z | \I \right) 
	\end{array}
\end{equation}
where $\int_{-\infty}^{\infty} d\mu_i$ denotes 
$\int_{-\infty}^{\infty} d\mu_i^1 \int_{-\infty}^{\infty} d\mu_i^2 ... 
\int_{-\infty}^{\infty} d\mu_i^p$ for components $1 .. p$ 
of $\mu_i$, with each integrated from $-\infty$ to $\infty$.  
The second line above is read as the probability of 
observing data X and clusters with means $\Mu$, given the 
cluster assignments $Z$, and $\I$. 
(This is very different to writing this as the product of probabilities of 
independent observations $X_i$, with cluster 
means $\mu_{Z_i}=\Pi_{g=1}^m \mu_g^{I(Z_i=g)}$.)
For independent normally distributed $\{ X_i\}$ with covariances 
$\{ \Gamma_i^{-1}\}$, sampled from clusters $g=1..m$ with means $\{ \mu_g\}$,
\begin{equation}\label{PPM}
	P\left( X | \Mu, Z, \I \right) 
	= \Pi_{g=1}^m \Pi_{i\in C_g} 
	f \left( x_i ; \mu_g, \Gamma_i \right)
\end{equation}
where $C_g=\{ i: Z_i=g\}$ are the members of the $g$th cluster, and,
\begin{equation}\label{toFact}
	f\left( x_i ; \mu_g, \Gamma_i \right)  
	=  
	\frac{1}{\sqrt{(2\pi)^p |\Gamma_i^{-1}|}}
	\exp \left( - \frac{1}{2} (x_i-\mu_g)^T 
	\Gamma_i (x_i-\mu_g) \right)
\end{equation}
Eq. \ref{PPM} has the form of a 
product partition model \cite{Hartigan1990,Barry1992}.  
In the following discussions we consider two priors for the 
means $\Mu$, one a flat prior with $P(\Mu|Z,\I)$ constant, 
and the other a normal distribution with, 
\begin{equation}\label{muPrior}
	P(\Mu|Z,\I) = \Pi_{g=1}^m f(\mu_g ; \mu_0, \Gamma_0 )
\end{equation}
We include Eq. \ref{muPrior} in the analysis below. 
Using Eqs. \ref{PPM} and \ref{muPrior}, 
Eq. \ref{bayes1} gives, 
\begin{equation}\label{LZ}
	\begin{array}{rl}
		P \left( Z | X, \I \right) 
		\propto 
		&  P(X | Z, \I ) P(Z| \I)
		\\
		\propto  
		&\int_{-\infty}^{\infty} d\mu_1 ... \int_{-\infty}^{\infty} d\mu_m
		\left(  \Pi_{i\in C_g}  f\left( x_i ; \mu_g, \Gamma_i \right)  \right)
		\left( \Pi_{g=1}^m  f(\mu_g ; \mu_0, \Gamma_0 ) \right)
		P(Z| I)
	\end{array}
\end{equation}
Now consider the likelihood for a {partition} \cite{Marcus} 
of type $(N_1, ... , N_m)$  
with $m$ clusters,  and $r_1 ... r_k$ clusters of size $1 ... k$. 
Notice that in Eq. \ref{LZ}, the factor $P(X|Z,\I)$ 
is invariant to permutations 
of cluster labels among clusters with the same number of items. 
In other words, provided $Z$ is from the same partition, 
then $P(X|Z,\I)$ is unchanged. 
Therefore we can write the probability of a partition $C$ as, 
\begin{equation}\label{clusterProb}
	\begin{array}{rl}
		P(C|X,\I)
		= \sum_{Z:Z\in C} P(Z|X,\I)
		&\propto \sum_{Z:Z\in C} P(X|Z,\I)P(Z|\I)
		\\
		&=P(X|Z\in C, \I)  \sum_{Z:Z\in C} P(Z|\I)
	\end{array}
\end{equation}
Where we used the above observation that 
$P(X|Z,\I)$ is the same for all $Z$ such that $Z\in C$.
Noting that the proposition $Z=\{Z_i=z_i\}$ is the same as 
$(Z \mbox{ and } M(Z)=m)$, where $M(Z)$ is the number of clusters, given 
the cluster memberships $Z$, then, 
\begin{equation}\label{PC0} 
	P( C| \I ) = 
	\sum_{Z:Z\in C} P(Z | \I) 
	= \sum_{Z:Z\in C} P( Z | M=m, \I ) P( M=m | \I ) 
\end{equation}
and assuming that all assignments $Z$ are equally likely given $M=m$, 
\begin{equation}\label{PCZ}
	\begin{array}{rl}
		\sum_{Z:Z\in C} P( Z | M=m, \I ) 
		&= \frac{\# \left(Z:Z\in C , M=m\right)}{\#\left(Z:M=m\right)}
		\\
		&=\frac{n!}{N_1! ... N_m! r_1! ... r_k!}\frac{1}{S(n,m)}
	\end{array}
\end{equation}
where $\#$ is used to denote ``the number of''. 
Eq. \ref{PCZ} equals the number of partitions of type $(N_1, ... , N_m)$, 
divided by the total number of partitions of $n$ identifiable 
items into $m$ non-empty clusters 
(the Stirling numbers of the second kind $S(n,m)$ \cite{Marcus}). 
An alternative calculation that leads to the same result is given in Appendix A.
For $P(M=m|\I)$, there are no labels $\{Z_i\}$ for items, or cluster 
numbers $\{N_i\}$.
Appendix B considers whether $P(M=m|\I)$ should reflect the 
number of ways of partitioning $n$ identical items into $m$ clusters, 
and concludes that this is unlikely to be an appropriate prior in most cases. 
Instead, the choice of $P(M=m|\I)$ might best be informed by the particular 
application.
For the examples here, we will take it as constant $P(M=m|\I)=1/n$.  
Continuing as before,
\begin{equation}\label{LZC}
	\begin{array}{rl}
		P \left( X | Z \in C, \I\right)
		&=\int_{-\infty}^{\infty} d\mu_1 ... \int_{-\infty}^{\infty} d\mu_m 
		P(X|\Mu,Z\in C, \I)P(\Mu|Z\in C,\I)
		\\
		&=
		\int_{-\infty}^{\infty} d\mu_1 ... \int_{-\infty}^{\infty} d\mu_m
		\left(  \Pi_{i\in C_g}  f\left( x_i ; \mu_g, \Gamma_i \right)  \right)
		\left( \Pi_{g=1}^m  f(\mu_g ; \mu_0, \Gamma_0 ) \right)
		\\
		&=
		\Pi_{g=1}^m 
		\int_{-\infty}^{\infty} d\mu_g f(\mu_g ; \mu_0, \Gamma_0 )  
		\left(  \Pi_{i\in C_g}  f\left( x_i ; \mu_g, \Gamma_i \right)  \right)
	\end{array}
\end{equation}
Putting Eqs. \ref{clusterProb}-\ref{PC} together we have, 
\begin{equation}\label{PC}
	\begin{array}{rl}
		P(C|X,\I) 
		\propto 
		& P(X|Z\in C, \I)P(C|\I)
		\\
		&=P(M=m|\I) \times 
		\frac{n!}{N_1! ... N_m! r_1! ... r_k!}
		\frac{1}{S(n,m)} 
		\\
		&\times \left( \Pi_{g=1}^m \Pi_{i\in C_g} \frac{1}{\sqrt{(2\pi)^p|\Gamma_i^{-1}|}} \right)\times
		\\    
		&\times \left( \Pi_{g=1}^m \int_{-\infty}^{\infty} d\mu_g 
		f_0(\mu_g) 
		\exp \left\{
		-\frac{1}{2}  \sum_{i\in C_g} 
		(x_i- \mu_g)^T \Gamma_i (x_i-\mu_g) \right\} \right) 
	\end{array}
\end{equation}
where $f_0(\mu_g)=f(\mu_g; \mu_0, \Gamma_0)$  is the prior for $\mu_g$.
The factors group into terms that correspond to the prior probability of  splitting $n$ equivalent items into between $1$ and $m$ non-empty clusters, 
multiplied by the prior probability of a partition of type $(N_1 ... N_m)$, 
given that there are $m$ clusters and making the prior assumption that 
all $Z$ are equally likely, multiplied by the probability of the data given 
that it has that partition. 
Eq. \ref{PC} has the reassuring quality that we could, quite reasonably, have written it 
down as the model we were going to study without any further justification. 
Using Eq. \ref{expansion}, the exponent in Eq. \ref{PC} may be rewritten as, 
\begin{equation}\label{expExp}
	\begin{array}{ll}
	\int_{-\infty}^{\infty} d\mu_g 
	f_0(\mu_g) 
	\exp \left\{
	-\frac{1}{2}  \sum_{i\in C_g} 
	(x_i- \mu_g)^T \Gamma_i (x_i-\mu_g) \right\} 
	&=
	\exp \left\{
	-\frac{1}{2}  \sum_{i\in C_g} 
	(x_i- \tilde{\mu}_g)^T \Gamma_i (x_i-\tilde{\mu}_g) 
	\right\} 
	\\
	&\times 
	\int_{-\infty}^{\infty} d\mu_g 
	f_0(\mu_g) 
	\exp \left\{
	-\frac{1}{2}  \sum_{i\in C_g} 
	(\mu_g- \tilde{\mu}_g)^T \tilde{\Gamma}_g (\mu_g-\tilde{\mu}_g) 
	\right\} 
	\end{array}
\end{equation}
where  
$\tilde{\Gamma}_g = \sum_{i\in C_g} \Gamma_i$, and 
$\tilde{\mu}_g = ( \sum_{i\in C_g} \Gamma_i )^{-1} \sum_{i\in C_g} \Gamma_i x_i$ 
as in Eqs. \ref{BmuG}-\ref{LambdaG}. 
When $f_0(\mu_g)$ is a flat prior, the integral over $\mu_g$ involving the last term   
gives $\sqrt{2\pi}^p \sqrt{|\tilde{\Gamma}_g^{-1}|}$ for each $g$. 
A normal prior with mean $\mu_0$ and covariance $\Gamma_0^{-1}$ has, 
\begin{equation}\label{normalPrior}
f_0(\mu_g) = \frac{1}{ \sqrt{2\pi}^p \sqrt{ |\Gamma_0^{-1}| } } 
	\exp 
	\left\{ 
		-\frac{1}{2} 
		\left( \mu_0 - \mu_g \right)^T 
		\Gamma_0 
		\left( \mu_0 - \mu_g \right)  
	\right\}
\end{equation}
Then for each $g=1 .. m$ in Eq. \ref{PC}, 
$f_0(\mu_g)$ leads to an extra factor of $1/\sqrt{2\pi}^p\sqrt{|\Gamma_0^{-1}|}$ 
and each sum over $C_g$ is modified to include $i=0$, 
with $x_0 \equiv \mu_0$, and $\Gamma_0^{-1}$ the covariance of the prior. 
For this case Eq. \ref{PC} becomes,
\begin{equation}\label{PofC2}
\begin{array}{rl}
P(C|X,\I) 
\propto 
& P(X|Z\in C, \I)P(C|\I)
\\
&=P(M=m|\I) \times 
\frac{n!}{N_1! ... N_m! r_1! ... r_k!}
\frac{1}{S(n,m)} 
\\
&\times \left( \Pi_{i=1}^n \frac{1}{\sqrt{(2\pi)^p|\Gamma_i^{-1}|}} \right)
\times \left( \Pi_{g=1}^m \sqrt{ \frac{\left| \tilde{\Gamma}^{-1}_g\right|}{\left| \Gamma_0^{-1} \right|} } \right) 
\\
&\times \Pi_{g=1}^m  
\exp \left\{
-\frac{1}{2}  \sum_{i\in C_g} 
(x_i- \tilde{\mu}_g)^T \Gamma_i (x_i-\tilde{\mu}_g) \right\} 
\end{array}
\end{equation}
The first three terms result from the prior $P(C| \I)$ for the 
partition based on combinatorial considerations. 
The fourth term on the right side of Eq. \ref{PofC2} is independent of the 
clustering model, the fifth term depends on the number and size of clusters and 
is discussed further in the next Section, and 
the sixth term is a sum of squares term that measures the 
goodness of fit of the data to the model.
If there is a flat prior, the extra factors of 
$1/\sqrt{2\pi}^p\sqrt{|\Gamma_0^{-1}|}$ are no longer present, 
leaving factors of $\sqrt{(2\pi)^p|\Gamma_g^{-1}|}$ 
in place of 
$\sqrt{ \left| \tilde{\Gamma}^{-1}_g\right|/\left| \Gamma_0^{-1} \right|}$. 
Compare Eq. \ref{PofC2} with the equivalent expression for 
a Gaussian mixture model with cluster means $\{\mu_g\}$ removed 
by marginalisation (integrated-out). 
Key differences are that: 
(i) a mixture model considers a distribution of named, identifiable clusters, 
not a partition, 
(ii) removing $\mu_g$ by marginalising  
the joint probability mass/density function in a Gaussian 
mixture model would lead to $n$ integrals 
over $\mu_g$, instead of just one for each cluster,  
(iii) here the covariances $\{\Gamma_i\}$ differ for each 
data point instead of only between clusters. 

\section{Bayesian information criteria}\label{bic} 

To understand the fifth term of Eq. \ref{PofC2} better, and to compare Eq. 
\ref{PofC2} with the Bayesian Information Criterion (BIC), 
use Eq. \ref{LambdaG} to write, 
\begin{equation}\label{ngp}
\left| \tilde{\Gamma}_g \right| 
= \left| n_g \frac{1}{n_g} \sum_{i\in C_g} \Gamma_i \right|
= \left| n_g \bar{\Gamma}_g \right|
= n_g^p \left| \bar{\Gamma}_g \right|
\end{equation}
with $\bar{\Gamma}$ now the mean of $\{ \Gamma_i \}$ for the cluster 
(including $\Gamma_0$ if there is a normally distributed prior). 
Using Eq. \ref{ngp} and $|A^{-1}|=1/|A|$, the fifth term of Eq. \ref{PofC2} 
can be written, 
\begin{equation}\label{meanI}
\Pi_{g=1}^m 
\sqrt{ \frac{| \tilde{\Gamma}^{-1}_g |}{| \Gamma_0^{-1} |} } 
 = \left( 
\Pi_{g=1}^m \frac{1}{\sqrt{n_g^{p}}} \right) 
\left( \Pi_{g=1}^m
\sqrt{ \frac{| \bar{\Gamma}^{-1}_g|}{| \Gamma_0^{-1} |} } \right)
\end{equation}
which is $O( \Pi_{g=1}^m 1/\sqrt{n_g^{p}} )$ for large $n_g$. 
Previous authors \cite{Jin1997, Chen1998} 
have argued for and used, an heuristic 
penalty term equivalent to $\log$ of Eq. \ref{meanI} 
when minimising the log-likelihood for clustering. 
In those studies the term was not derived, but was one of several 
proposed alternatives \cite{Jin1997, Chen1998}.  
Substituting Eq. \ref{meanI} into \ref{PofC2}, then taking the logarithm  
and multiplying by $-2$, gives, 
\begin{equation}\label{logLik}
\begin{array}{ll}
-2 \log\left( P\left( C | X, \Gamma  \right)  \right) 
=&-2\log \left( P(M=m|\I) \right)
\\&
-2\log\left( \frac{n!}{N_1! ... N_m! r_1! ... r_k!} \right)
\\&
+2\log\left( S(n,m) \right) 
\\&
+ \sum_{i=1}^n \log  \left| 2\pi \Gamma_i^{-1} \right| 
\\&+ 
\sum_{g=1}^m
\sum_{i \in C_g} 
(x_i- \tilde{\mu}_g)^T \Gamma_i (x_i-\tilde{\mu}_g) 
\\&+ 
\sum_{g=1}^m  p \log n_g
\\&
- \sum_{g=1}^m 
\log \left( |\bar{\Gamma}_g^{-1}|/|\Gamma_0^{-1}| \right)
\end{array} 
\end{equation}
The maximum likelihood estimates for $\{\mu_g\}$ can be found from 
Eq. \ref{PC} by taking derivatives with respect to the components 
of each $\mu_g$, and have $\hat{\mu}_g=\tilde{\mu}_g$. 
Hence the fourth 
and fifth terms are the sum of $-2$ times the maximum-likelihood estimates for 
the log-likelihoods\footnote{Strictly speaking, we are discussing the 
	log of the posterior distribution and the 
	``MLE'' is the maximum a posteriori probability (``MAP'') estimate.} 
of each cluster ($-2\log \hat{L}$), 
the first three terms arise from combinatorial considerations, 
and the final two terms asymptote to 
$\sum_{g=1}^m p \log n_g$ as $n_g \rightarrow \infty$. 
For the limit of large $n_g$ this becomes,
\begin{equation}\label{clusterBIC}
-2 \log\left( P\left( C | X, \Gamma  \right)  \right)  
= -2\log P(C|\I) + \sum_{g=1}^m \left\{ - 2\log(\hat{L}_g) + p \log n_g \right\}
+ O \left( 1 \right)
\end{equation}
where the terms involving 
$\log \left( |\bar{\Gamma}_g^{-1}|/|\Gamma_0^{-1}| \right) \sim 1$, and 
are not explicitly included to allow comparison with the 
Bayesian Information Criterion (BIC). 
Usually $|\Gamma_0^{-1}|>|\bar{\Gamma}_g^{-1}|$ and the terms neglected 
from Eq. \ref{logLik} would be positive, the same sign as $p\log(n_g)$ and 
will penalise larger numbers of clusters,  
but their size would depend 
on the details of the clustering. 
Neglecting these terms will usually underestimate the penalty 
associated with having more clusters. 
If there is a normal prior, then $n_g$ equals one plus the number of 
elements in a cluster,  and $\log(n_g)\rightarrow 2$ as $m \rightarrow n$. 
The BIC for a parametric model is usually defined as \cite{Wasserman2010}, 
\begin{equation}\label{BIC0}
BIC = -2 \log \left( \hat{L} \right) + k \log (n)
\end{equation}
where $\hat{L}$ is the log-likelihood at the MLE, $k$ is the number of parameters, 
and $n$ is the number of data. 
If we were simply studying a single cluster with a flat prior, then 
Eqs. \ref{clusterBIC} and \ref{BIC0} would be identical, but with $k=p$. 
If $n_g \gg 1$ for all clusters, then Eqs. \ref{logLik} and \ref{clusterBIC}, 
are equivalent to the sum of BICs for each cluster, plus a combinatorial term for 
the probability of sampling the set of clusters by chance (that requires a prior for the number of clusters $P(M=m|\I)$). 
Note that in general $n_g$ is not always large, and in many applications 
a cluster can represent a single item. 
In clustering applications the usual BIC approximation is not applicable unless 
all clusters are large with $n_g \gg 1$, but when the data are normally 
distributed, then the exact posterior distribution can be used 
(Eqs. \ref{PofC2} and \ref{logLik}). 

\section{Examples}

A recent epidemiological study using UK Biobank data \cite{Webster2021}, 
estimated associations with 12 well-known risk factors in 
over $400$ diseases using a proportional hazards survival analysis. 
Full details of the study and dataset are in Ref. \cite{Webster2021}, 
summary statistics and code are at: osf.io.
Diseases with statistically significant differences between men and women after 
an FDR multiple-testing adjustment \cite{Wasserman2010} were excluded, as were 
diseases that failed a global $\chi^2$ test of the proportional hazards model 
using Schoenfeld residuals, and only diseases whose associations remained 
statistically significant after a Bonferroni adjustment were kept. 
This left  78 pairs of diseases affecting men and women (156 diseases in total).

\begin{figure}[ht]
	\centering
	\includegraphics[width=0.85\linewidth,angle=0]{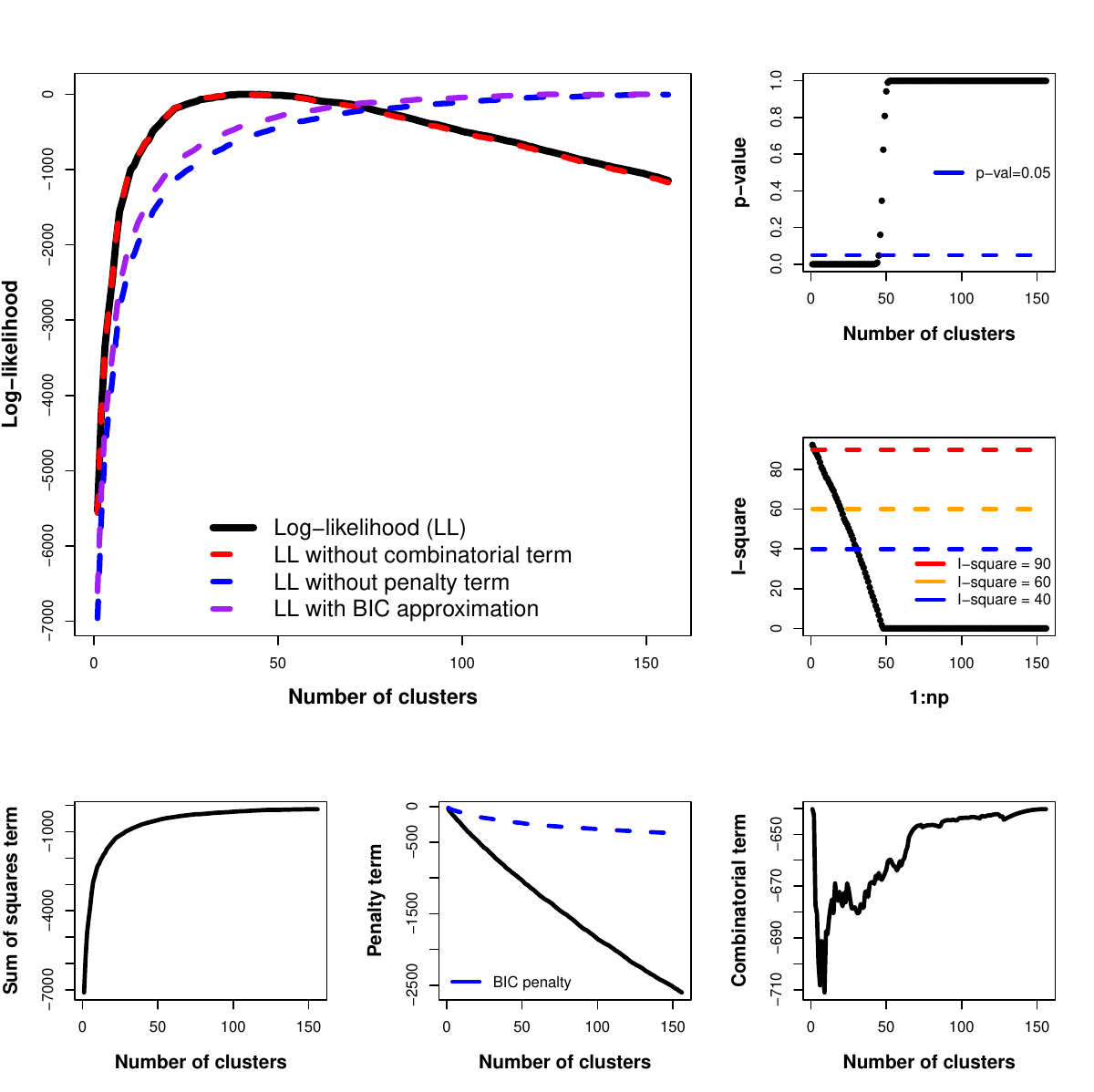}
	\caption[]{
		Using Eq. \ref{PofC2}, log-likelihoods are calculated and their 
		maxima subtracted (top left): 
		all terms (black), without the combinatorial term $\log P(C|\I)$, without   
		the penalty term 
		$\log \Pi_{g=1}^m \sqrt{|\tilde{\Gamma}^{-1}|/|\Gamma_0^{-1}|}$, 
		and with the penalty term replaced by the usual BIC approximation 
		$-\sum_{g=1}^m p \log n_g$.  
		Individual terms are plotted 
		in the bottom three figures, including the sum of squares term 
		$-(1/2)\sum_{g=1}^m
		\sum_{i \in C_g} 
		(x_i- \tilde{\mu}_g)^T \Gamma_i (x_i-\tilde{\mu}_g)$. 
		The p-value for heterogeneity and the $I^2$ statistic 
		Eqs \ref{clusteringQsq} and \ref{Isq}, are the top-right figures. 
		The data describe associations between common risk factors and 
		156 diseases \cite{Webster2021}, hierarchically clustered to 
		to search for shared aetiologies. 
	}\label{figLL}
\end{figure}
The authors wished to cluster diseases using MLEs for associations between 
exposures and disease incidence, because exposure-disease associations 
were expected to reflect causal disease pathways.  
Because associations between diseases and risk factors are often strongly 
correlated, the covariances of MLEs must be accounted for when clustering. 
Due to very different incident rates (sample sizes), covariances were expected 
to differ substantially even if diseases originated from the same cluster. 
Therefore the authors used the Bhattacharyya distance and 
hierarchical clustering, the latter allowing easier interpretation of the 
resulting clusters. 
A limitation of the approach, is that it did not determine how many 
clusters to consider.
The authors \cite{Webster2021} used the ad-hoc elbow criterion to keep 24 clusters, 
but acknowledged the need for an objective selection criterion. 
The problem was to cluster $\{\hat{\mu}_i, \hat{\Sigma}_i \}$ into groups 
with similar $\{\mu_i\}$. 
A simple model that accounts for the uncertainty of estimates 
$\{\hat{\mu}_i\}$ through their covariances, is to take $\hat{\mu}_i \sim N(\mu_g,\hat{\Sigma}_i)$ for clusters with means $\{\mu_g\}$; which is the 
model considered here. 
Figure \ref{figLL} shows the log-likelihood for hierarchical clustering with the Bhattacharyya distance (see Appendix D), and a normal prior with an isotropic 
covariance $\Sigma_0$ with $\sqrt{tr(\Sigma_0)}=3$  
(when $\mu_g\sim N(\mu_0,\Sigma_0)$, $E[(\mu_g-\mu_0)^2]=tr(\Sigma_0)$). 
The choice of $\sqrt{tr(\Sigma_0)}=3$ was intended to be large enough to 
included the largest expected relative risks from smokers and lung-cancer, 
that has been found to be \cite{Pirie2013} of order $20 \simeq \exp(3)$. 
The figure plots the log-likelihood as calculated using Eq. \ref{logLik}, whose 
maximum is at $42$ clusters. 
It also shows the influence of the combinatorial term $\log P(C|\I)$ (which 
was small), the penalty term 
$\log \Pi_{g=1}^m \sqrt{|\tilde{\Gamma}^{-1}|/|\Gamma_0^{-1}|}$, 
and the log-likelihood using the equivalent term from the usual 
BIC approximation 
$-(1/2)\sum_{g=1}^m p \log n_g$.  
As expected, as the number of clusters increase and the typical cluster 
sizes become smaller, the BIC approximation becomes worse.
The p-value for heterogeneity (Eq. \ref{clusteringQsq}), rises above $0.05$ 
at $46$ clusters. 
An $I^2$ value less than $40$ is usually regarded as a low level of 
heterogeneity \cite{Borenstein2009}, 
and this occurs for $31$ or more clusters, 
consistent with the maximum log-likelihood at $42$. 
With $42$ clusters, $72$\% of the $78$ disease types (78 in men and 78 in 
women), appeared with their opposite-sex pair in 
the same cluster.
The diseases and their clusters are listed in the Supplementary 
Material, along with the results from a sensitivity analysis 
that is described in Section \ref{CIs}. 
Further examples using simulated data are in the Supplementary Material. 
These include sensitivity analyses that are described in 
Section \ref{CIs}, some explorations for how results are influenced by 
noise in the data, and specific numerical test cases.

\section{Discussion}\label{discussion}

\subsection{Limits of a normal prior}\label{limits}

To better understand the influence of the normal prior $f_0(\mu_g)$, take 
$\Gamma_0^{-1}=I \sigma_0^2$ where $I$ is the identity matrix, and explore 
the limits of a sharply peaked prior with $\sigma_0^2 \rightarrow 0$ and 
of a flat prior with $\sigma_0^2 \rightarrow \infty$. 
These limits concern the last three terms in Eq. \ref{logLik}. 
Firstly consider $\sigma_0^2 \rightarrow \infty$, a limit that emphasises that 
a flat prior cannot be considered as a limiting case of a normal prior.
In the final term 
$\tilde{\Gamma}_g = \sum_{k\in C_g} \Gamma_k \rightarrow  \sum_{k\in C_g, k\neq 0} \Gamma_k$. 
In the fifth term the components of $\Gamma_0 = I/\sigma_0^2 \rightarrow 0$, and, 
\begin{equation}
\tilde{\mu}_g \rightarrow 
\left( \sum_{i\in C_g, i\neq 0} \Gamma_i \right)^{-1} 
\sum_{i\in C_g, i\neq 0} \Gamma_i x_i
\end{equation}
and $-\log \left| \Gamma_0^{-1} \right|$ 
becomes $-(p/2)\log(2\pi \sigma_0^2)$ where $p$ is the 
dimension of $\mu$, and $\log(\sigma_0^2) \rightarrow \infty$ 
as $\sigma_0^2 \rightarrow \infty$. 
The divergent behaviour arises from requiring that $f_0(\mu_g)$ is 
correctly normalised, and can be understood from Eqs. \ref{PC} and \ref{normalPrior}. 
As $\sigma_0^2$ becomes larger, the maximum of $f_0(\mu_g)$ must become 
increasingly small to ensure that it is correctly normalised, and there is an 
extra factor for each cluster. 
For this example with $\Gamma_0^{-1}=\sigma_0^2 I$, the first term provides a penalty that 
is proportional to the number of free parameters $m\times p$. 
For the alternative limit with $\sigma_0^2 \rightarrow 0$, Appendix \ref{sharpPrior} 
shows that the fifth and seventh terms can be combined and will tend to zero, and 
the sixth term tends to, 
\begin{equation}
- \frac{1}{2} 
\sum_{g=1}^m
\sum_{i \in C_g} 
(x_i- \tilde{\mu}_g)^T \Gamma_i (x_i-\tilde{\mu}_g) 
= - \frac{1}{2} \sum_{g=1}^m
\sum_{i \in C_g, i\neq 0} 
(x_i- {\mu}_0)^T \Gamma_i (x_i-{\mu}_0) + O(\sigma_0^2)
\end{equation}
so that the log-likelihood is equivalent to clusters that all have 
the same centres $\mu_g=\mu_0$, as we might have expected.  

\subsection{Equal covariances, and relation to k-means}\label{kmeans}

Consider equal covariances for all the data 
points\footnote{For example, this might occur for longitudinal epidemiological 
	data with repeated measurements to allow estimation of either measurement errors 
	or of the intrinsic variation within individuals. 
	An additional set of measurements at a later time allows a covariance to 
	be calculated, and these can be averaged across all individuals to give 
	an estimated covariance for within-person measurements. 
	The average ``within-person'' covariance is different to 
	the ``population'' covariance of data for a single time point, for example because the 
	variation within an individual could be less than between individuals.}, 
with $\Gamma_i^{-1}=\Gamma^{-1}$ and a flat prior.
Using $|AB|=|A||B|$ and $|B^{-1}|=1/|B|$,
the final term in Eq. \ref{logLik} simplifies to, 
\begin{equation}
- \frac{1}{2} \sum_g 
\left( \log \left|\bar{\Gamma}_g^{-1} \right| \right) 
=
\frac{1}{2} \sum_g  \log \left| \Gamma \right| 
= 
\frac{m}{2} \log \left| \Gamma \right| 
\end{equation}
With $\Gamma_i=\Gamma$ and $\tilde{\Gamma}= n_g \Gamma$, then 
$\tilde{\mu}=\frac{1}{n_g}\sum_{i\in C_g} x_i$, and 
the fifth term in Eq. \ref{logLik} is,
\begin{equation}\label{kmss}
\sum_{g=1}^m \sum_{i\in C_g} (x_i -\tilde{\mu}_g)^T \Gamma_i (x_i -\tilde{\mu}_g) 
=
\sum_{g=1}^m \sum_{i\in C_g} 
\left(x_i - \frac{1}{n_g} \sum_{j\in C_g} x_j\right)^T 
\Gamma 
\left(x_i - \frac{1}{n_g} \sum_{j\in C_g} x_j\right)
\end{equation}
Eq. \ref{expansion3} shows how Eq. \ref{kmss} relates to the sum of 
pairwise within-cluster differences. 
The log-likelihood is, 
\begin{equation}\label{logLikKmeans}
\begin{array}{ll}
-2 \log\left( P\left( C | X, \Gamma  \right)  \right) 
=&-2\log \left( P(C|\I) \right)
\\&
+ \sum_{i=1}^n \log  \left| 2\pi \Gamma^{-1} \right| 
\\&+ 
\sum_{g=1}^m \sum_{i\in C_g} 
\left(x_i - \frac{1}{n_g} \sum_{j\in C_g} x_j\right)^T 
\Gamma 
\left(x_i - \frac{1}{n_g} \sum_{j\in C_g} x_j\right)
\\&+ 
\sum_{g=1}^m  p \log n_g
\\&
+ m \log  |\Gamma| 
\end{array} 
\end{equation}
For a flat prior and equal covariances, maximising the log-likelihood is 
similar to minimising the sum of squares as for k-means, but 
with two extra terms to penalise the model's complexity. 
The penalty terms for $-2\times$ the log-likelihood are 
$-2\log \left( P(C|\I) \right)$, and, 
\begin{equation}
\sum_{g=1}^m  \left( p\log n_g + \log |\Gamma | \right)
\end{equation}
The second term's sign depends on whether $0<|\Gamma|<1$ or $1<|\Gamma|$, and 
can penalise more, or fewer clusters, with large covariances ($0<|\Gamma|<1$) 
favouring more clusters but smaller covariances ($1<|\Gamma|$) favouring 
fewer clusters. 
Comparison with Eq. \ref{PofC2} indicates that $|\Gamma|>1$ 
(that favours fewer clusters), 
is analogous having a prior with covariance $\Sigma_0$ that 
has $|\Sigma|<|\Sigma_0|$ (with $|\Gamma| > |\Gamma_0|$), as would usually be the case.  
The first term is an entropy-like term, and the concave shape of $\log n_g$ 
will penalise similarly sized clusters (when summed over $g$), 
with $ \frac{1}{m} \sum_{g=1}^m \log n_g \leq \log \frac{1}{m} \sum_{g=1}^m n_g 
= \log(n/m)$.  
Unfortunately the k-means algorithm will not minimise Eq. \ref{logLikKmeans} 
by assigning items to clusters with the nearest mean, 
because the penalty terms depend on the number and composition of clusters.

\subsection{Heterogeneity, composite endpoints, and sub-group analysis}

Composite endpoints consist of several grouped symptoms or diseases, and 
are intrinsic to how diseases are defined and studied. 
Since the first statistical studies of disease by John Graunt in the 
1600s \cite{Graunt}, there has 
been a trade-off between definitions that are sufficiently specific to 
distinguish different underlying disease processes, and sufficiently broad to 
allow a meaningful statistical study. 
This is particularly apparent in clinical trials and epidemiological studies 
where data are costly or unavailable. 
Large population datasets with detailed genetic and biological information 
are providing new data-driven definitions of disease, identifying distinct 
subtypes of disease, and collections of diseases with potentially shared 
underlying causes \cite{Webster2021}. 
Statistical methods can assess whether a composite endpoint consisting of several 
potentially distinct diseases, is consistent with 
its assumed properties, such as testing the constituent diseases for heterogeneity 
of their disease-risk associations. 
Heterogeneity is conventionally tested with a $Q$ or $I^2$ statistic \cite{Borenstein2009}. 
When heterogeneity is anticipated in advance, potential subgroups are 
often proposed as an alternative to a single cluster. 
Because the subgroups are pre-specified, as opposed to clustered, the larger 
number of free parameters does not guarantee a less heterogeneous result. 
In principle Eq. \ref{logLik} can, and should, be used in preference to the 
Q-statistic when deciding whether a proposed sub-group is a better representation 
of the data than a single cluster, or a different sub-group. 
However, because the sub-groups are pre-specified in advance, the first 
three (combinatoric terms) should not be included.  
Clustering algorithms usually ensure that 
the $Q$ or $I^2$ statistics decrease as the 
number of clusters are increased. 
In contrast, because Eq. \ref{logLik} accounts for the number of free parameters, 
it can have have a minimum at a particular number of clusters. 
This suggests an alternative approach that does not require pre-specified 
subgroups, but instead tests if Eq. \ref{logLik} is minimised by one, or more 
clusters. 
The merits of this approach for applications such as meta-analyses, will need 
exploring in greater detail elsewhere. 

\subsection{Statistical tests and confidence sets}\label{CIs}

Maximum likelihood estimates are usually reported with a confidence set 
to provide a measure of uncertainty in the estimate. 
In principle this is possible for clusterings, either with confidence sets for 
the overall configuration of clusters, or for the composition of individual clusters. 
One difficulty is that Eq. \ref{logLik} does not have the usual asymptotic properties of 
log-likelihoods near the maximum likelihood estimate (MLE), that would usually involve a sum over independent random variables, for which the score function asymptotically  has  a normal distribution. 
In contrast, Eq. \ref{logLik} involves a sum over clusters and the cluster's membership. 
In principle, the Bootstrap method offers a simple way to generate confidence sets.
Data can be randomly sampled with replacement as usual, and an optimal clustering found 
by a suitable method. 
However, a problem with this simple approach, is that the number of clusters can change due 
to some clusters not being sampled.
The underlying issue is that when bootstrap is usually used, each data point contributes  
information about all the parameters, e.g. the intercept and slope of a line. 
When clustering, a data point only contributes information about the 
cluster it belongs to.
It might be possible to work around these issues with a more complex sampling method, 
or a more careful interpretation of results. 
An alternative option when the data and sufficient computing power are available,  
is to bootstrap sample the underlying data used to generate the MLEs that 
are subsequently clustered, to obtain the optimum clustering and 
log-likelihood for each sample. 
This will give a log-likelihood $l_i$, and other properties 
for each clustered sample, such as the number of clusters $m_i$. 
After generating sufficient samples, a confidence set can then be formed  
using the empirical distribution of the log-likelihoods $\{l_i\}$, and the  
properties of samples within the confidence set can be studied. 
For example, the distribution for the number of clusters can then  
reported with a confidence set. 
The ``Jacknife'' \cite{Wasserman2010} can provide a simple estimate for variances. 
In practice it is often a poor estimate when a statistic is not smooth, 
as will be the case for statistics such as the MLE for the number 
of clusters. 
However, the variation in statistics under the leave-one-out procedure can 
provide a simple but useful indication of how sensitive our MLE clustering 
estimate is to small changes in the data. 
Although the Jacknife procedure is unsuitable for estimating a confidence interval 
for the number of clusters, the leave-one-out procedure can provide a valuable 
sensitivity analysis, that indicates when there is uncertainty in the MLE for the 
number of clusters. 
A histogram formed by systematically removing one item 
at a time before re-clustering using the Bhattacharyya distance, 
hierarchical clustering, and Eq. \ref{PofC2}, 
is discussed further in the Supplementary Material. 
The majority of cases had between $39$ and $43$ clusters, with a 
peak at $42$. 
The Supplementary Material includes a more systematic study 
using simulated data, that explores how the histograms change as noise 
in the data is gradually increased. 

\subsection{Model improvements - uncertainty in covariance estimates}\label{covEst}

The clustering model does not account for uncertainty in the 
estimated covariances, that for MLEs, are estimated from the underlying 
data. 
One model to account for uncertainty, is to model the 
estimated covariances $\hat{\Sigma}_i$ as sampled from a Wishart distribution  
$W(\hat{\Sigma}_i | \Sigma_i/n_i,n_i)$, where $n_i$ are the 
number of data in the estimate for $\hat{\Sigma}_i$, and $\Sigma_i$ is 
the unknown covariance. 
Then using Bayes theorem, 
\begin{equation}
P\left( \Sigma_i | \hat{\Sigma}_i, n_i \right) 
P\left(  \hat{\Sigma}_i \right) = 
W \left( \hat{\Sigma}_i | \Sigma_i/n_i, n_i \right) 
P\left( \Sigma_i \right)
\end{equation}
With Jeffrey's prior $P(\Sigma_i) \propto |\Sigma_i|^{{-p+1}/2}$, it can be 
shown that \cite{Sellentin2015}, 
\begin{equation}\label{Tdist}
P(X_i | \mu, \hat{\Sigma}_i, n_i) = \frac{c_p}{\sqrt{|\hat{\Sigma_i}|}}
\frac{1}{\left[ 1 + \frac{(X_i-\mu)^T\hat{\Sigma}_i^{-1}(X_i-\mu)}{n_i} \right]^{\frac{n_i+1}{2}}}
\end{equation}
where,
\begin{equation}
c_p = \frac{\Gamma((n_i+1)/2)}{[\pi n_i]^{p/2} \Gamma(\frac{n_i+1-p}{2})}
\end{equation}
Eq. \ref{Tdist} is a form of multivariate t-distribution, and can be used to 
calculate the likelihood. 
The author is presently unaware of a suitable conjugate prior or 
generalisation of the log-likelihood calculation of the previous Sections. 
An option is to formulate a numerical calculation or to implement a 
Dirichlet Process Mixture Model (DPMM) \cite{Ross2020}.

\section*{Conclusions}

This research arose from a need to combine the best existing epidemiological 
methods with clustering techniques, so as to identify shared causes of 
disease \cite{Webster2021}.  
This was accomplished by using established parametric survival models 
to characterise the data, through MLEs for associations between exposures 
and disease risks. 
MLEs are (asymptotically) normally distributed, which led to 
the general problem of clustering (multivariate) normally distributed data, to  
determine the number and composition of clusters. 
The posterior distribution for this model was calculated by marginalising the 
unknown cluster centres to give Eqs. \ref{PofC2} and \ref{logLik}, a procedure 
that is usually combined with a Laplace approximation when calculating the BIC. 
In the limit where the number of items in each individual cluster is large 
enough, then Eq. \ref{logLik}, will asymptotically agree with an 
equivalent expression using a sum of the usual BIC estimated for each 
cluster, plus a combinatoric term. 
The combinatoric term is intended to indicate the probability of the clusters 
occurring by chance. 
In general it is inappropriate to use the usual BIC to compare clusterings, 
and when used with a normal prior it will usually underestimate 
the penalty associated with having more clusters, but when the underlying 
data are normally distributed MLEs then -
Eqs. \ref{PofC2} and \ref{logLik} can be used. 

\appendix 

\section{Multivariate (inverse variance weighted) sums of squares I}\label{mivw}

Note that the covariances and their inverses are symmetric, and expand,
\begin{equation}\label{expansion}
\begin{array}{l}
\sum_i \left( x_i - \mu_g \right)^T \Gamma_i \left( x_i -\mu_g \right)
\\= 
\sum_i 
\left( \left( x_i - \tilde{\mu}_g \right) - 
\left( \mu_g - \tilde{\mu}_g \right) \right)^T 
\Gamma_i 
\left( \left( x_i - \tilde{\mu}_g \right) - 
\left( \mu_g - \tilde{\mu}_g \right) \right)
\\=
\sum_i 
\left( x_i - \tilde{\mu}_g \right)^T 
\Gamma_i 
\left( x_i - \tilde{\mu}_g \right)
+
\sum_i 
\left( \mu_g - \tilde{\mu}_g \right)^T 
\Gamma_i 
\left( \mu_g - \tilde{\mu}_g \right)
\\-2
\sum_i 
\left( \mu_g - \tilde{\mu}_g \right)^T 
\Gamma_i 
\left(  x_i - \tilde{\mu}_g  \right)
\\=
\sum_i 
\left( x_i - \tilde{\mu}_g \right)^T 
\Gamma_i 
\left( x_i - \tilde{\mu}_g \right)
+
\left( \mu_g - \tilde{\mu}_g \right)^T 
\left( \sum_i \Gamma_i \right)
\left( \mu_g - \tilde{\mu}_g \right)
\\-2
\left( \mu_g - \tilde{\mu}_g \right)^T 
\left( 
\left( \sum_i  \Gamma_i  x_i \right) 
- 
\left( \sum_i \Gamma_i \right) \tilde{\mu}_g 
\right)
\end{array}
\end{equation}
where the sums are over all $i$ in cluster $C_g$. 
If $\tilde{\mu}$ takes the specific form given by Eq. \ref{BmuG}, with, 
\begin{equation}\label{tildeMu}
		\tilde{\mu}_g =   
	\left( \sum_{i \in C_g} \Gamma_i \right)^{-1} 
	\left( \sum_{i \in C_g} \Gamma_i x_i \right)
\end{equation}
then the terms $\left( \sum_i  \Gamma_i  x_i \right)$ and 
$\left( \sum_i \Gamma_i \right) \tilde{\mu}_g$ in the last term of the 
final line  cancel.
Writing $\tilde{\Gamma}_g$ as in Eq. \ref{LambdaG}, with, 
\begin{equation}
	\tilde{\Gamma}_g = \sum_{i \in C_g} \Gamma_i
\end{equation}	
then the resulting equation becomes, 
\begin{equation}\label{expansion2}
	\begin{array}{l}
		\sum_{i \in C_g} \left( x_i - \mu_g \right)^T \Gamma_i \left( x_i -\mu_g \right)
		\\= 
		\sum_{i \in C_g} 
		\left( x_i - \tilde{\mu}_g \right)^T 
		\Gamma_i 
		\left( x_i - \tilde{\mu}_g \right)
		+
		\left( \mu_g - \tilde{\mu}_g \right)^T 
		\tilde{\Gamma}_g 
		\left( \mu_g - \tilde{\mu}_g \right)
	\end{array}
\end{equation}
that can be rearranged to give Eq. \ref{penult}. 

\section{Multivariate sums of squares II}\label{sumOfSquares}

To integrate over $\mu$ write $\Gamma_i=\Sigma_i^{-1}$ and note 
that $\Sigma_i$ and their inverses $\Gamma_i$ are symmetric, and use this to write, 
\begin{equation}\label{sseq}
	\begin{array}{l}
		\sum_i (x_i-\mu)^T  \Gamma_i (x_i-\mu)
		\\ = 
		\left( \mu - \left( \sum_i \Gamma_i \right)^{-1} \sum_i \Gamma_i x_i \right)^T 
		\left( \sum_i \Gamma_i \right) 
		\left( \mu - \left( \sum_i \Gamma_i \right)^{-1} \sum_i \Gamma_i x_i \right)
		\\+ \sum_i x_i^T \Gamma_i x_i 
		- \left( \sum_i x_i^T \Gamma_i \right) 
		\left( \sum_i \Gamma_i \right)^{-1} 
		\left( \sum_i \Gamma_i x_i \right)
	\end{array}
\end{equation}
The terms involving $\mu$ in \ref{sseq} will factorise in 
Eq. \ref{PC}, and lead to Gaussian integrals that can be 
integrated to give functions of $\{\Gamma_i\}$ 
that are independent of $\{ x_i \}$. 
The remaining terms are, 
\begin{equation}\label{remEq}
	\begin{array}{l} 
		\sum_i x_i^T \Gamma_i x_i 
		- \left( \sum_i x_i^T \Gamma_i \right)
		\left( \sum_k \Gamma_k \right)^{-1}
		\left( \sum_j \Gamma_j x_j \right) 
		\\=
		\sum_{i,j} x_i^T 
		\Gamma_i \left( \sum_k \Gamma_k \right)^{-1} \Gamma_j 
		x_i
		- \sum_{i,j} x_i^T 
		\Gamma_i \left( \sum_k \Gamma_k \right)^{-1} \Gamma_j
		x_j
	\end{array}
\end{equation}
Because $\Sigma_i$ and their inverses $\Gamma_i$ are symmetric, then $C_{ij}=\Gamma_i \left( \sum_k \Gamma_k \right)^{-1}\Gamma_j$ 
has $C_{ij}=C_{ji}^T$, as can be seen by taking the transpose of $C_{ij}$. 
Using $C_{ij}=C_{ji}^T$, $a^T b=b^T a$ for vectors $a$ and $b$,  and relabeling the indices $i$ and $j$, 
\begin{equation}\label{Cij}
	\begin{array}{l} 
		\frac{1}{2} \sum_{i,j} (x_i-x_j)^T C_{ij} (x_i-x_j) 
		\\= 
		\frac{1}{2} \sum_{i,j} x_i^T C_{ij} x_i 
		+\frac{1}{2} \sum_{i,j} x_j^T C_{ij} x_j 
		-\frac{1}{2} \sum_{i,j} x_i^T C_{ij} x_j 
		-\frac{1}{2} \sum_{i,j} x_j^T C_{ij} x_i
		\\= 
		\sum_{i,j} x_i^T C_{ij} x_i 
		-\frac{1}{2} \sum_{i,j} x_i^T C_{ij} x_j 
		-\frac{1}{2} \sum_{i,j} x_j^T C_{ji}^T x_i
		\\=
		\sum_{i,j} x_i^T C_{ij} x_i - x_i^T C_{ij} x_j
	\end{array}
\end{equation}
where the last line is the same form as the last line of Eq. \ref{remEq}. 
Hence using Eqs. \ref{sseq}, \ref{remEq}, and \ref{Cij}, we have, 
\begin{equation}\label{ssEqf2}
	\begin{array}{l} 
		\sum_i (x_i-\mu)^T  \Gamma_i (x_i-\mu)
		\\ = 
		\left( \mu - \left( \sum_i \Gamma_i \right)^{-1} \sum_i \Gamma_i x_i \right)^T 
		\left( \sum_i \Gamma_i \right) 
		\left( \mu - \left( \sum_i \Gamma_i \right)^{-1} \sum_i \Gamma_i x_i \right)
		\\+
		\frac{1}{2} \sum_{i,j} (x_i-x_j)^T 
		\Gamma_i \left( \sum_k \Gamma_k \right)^{-1}  \Gamma_j
		(x_i-x_j) 
	\end{array}
\end{equation}
where the sums over $i$, $j$, and $k$ will range over elements in cluster $g$. 
Using Eqs. \ref{BmuG} and \ref{LambdaG}, Eq. \ref{ssEqf2} can alternately be written as,
\begin{equation}\label{ss2}
	\begin{array}{l} 
		\sum_{i \in C_g} (x_i-\mu)^T  \Gamma_i (x_i-\mu)
		\\ = 
		\left( \mu - \tilde{\mu}_g \right)^T 
		\tilde{\Gamma}_g 
		\left( \mu - \tilde{\mu}_g \right)
		\\+
		\frac{1}{2} \sum_{i,j \in C_g} (x_i-x_j)^T 
		\Gamma_i \tilde{\Gamma}_g^{-1} \Gamma_j
		(x_i-x_j) 
	\end{array}
\end{equation}
and comparing with Eq. \ref{expansion2}, we can infer that,
 \begin{equation}\label{expansion3}
 	\begin{array}{l} 
		\sum_{i \in C_g} 
		\left( x_i - \tilde{\mu}_g \right)^T 
		\Gamma_i 
		\left( x_i - \tilde{\mu}_g \right)
		=
 		\frac{1}{2} \sum_{i,j \in C_g} (x_i-x_j)^T 
 		\Gamma_i \tilde{\Gamma}_g^{-1} \Gamma_j
 		(x_i-x_j) 
 	\end{array}
 \end{equation}
as can be confirmed by expanding out and using the definitions 
of $\tilde{\mu}$ and $\tilde{\Gamma}_g$ given by Eqs. \ref{BmuG} and \ref{LambdaG}.

\section{Using Bayes theorem to estimate the cluster mean}\label{BayesApp}

Bayes theorem gives, 
\begin{equation}
P(\mu_g | X ) = \frac{ P(X|\mu_g) P(\mu_g)}{P(X)}
\end{equation}
Because $\int P(\mu_g | X) d\mu_g =1$, this may alternately be written as, 
\begin{equation}\label{b1}
P(\mu_g | X ) = \frac{ P(X|\mu_g) P(\mu_g)}{\int d\mu_g P(X|\mu_g) P(\mu_g)}
\end{equation}
Using Eq. \ref{expansion2} we can write, 
\begin{equation}
\begin{array}{ll}
P(X|\mu_g) P( \mu_g) 
& \propto \exp \left\{ -\frac{1}{2} 
\sum_{i \in C_g} \left( x_i - \mu_g \right)^T \Gamma_i \left( x_i -\mu_g \right)
\right\}
\\
&= \exp \left\{ -\frac{1}{2} 
\sum_{i \in C_g} \left( x_i - \tilde{\mu}_g \right)^T \Gamma_i \left( x_i -\tilde{\mu}_g \right)
\right\}
\\
& \times \exp \left\{ -\frac{1}{2} 
\left( \mu_g - \tilde{\mu}_g \right)^T \tilde{\Gamma}_g \left( \mu_g -\tilde{\mu}_g \right)
\right\}
\end{array}
\end{equation} 
The factors involving $x_i$ are independent of $\mu_g$, and will cancel each 
other in the 
numerator and denominator of Eq. \ref{b1}. 
Integrating over $\mu_g$ in the denominator then leads to, 
\begin{equation}
P(\mu_g|X)= 
\frac{\exp\left\{ 
	- \frac{1}{2} 
	\left( \mu_g - \tilde{\mu}_g \right)^T \tilde{\Gamma}_g \left( \mu_g -\tilde{\mu}_g \right)\right\}}{\sqrt{2\pi}^p \sqrt{|\tilde{\Gamma}_g^{-1}|}}
\end{equation}
where $p$ is the dimension of $\mu_g$. 

\section{Alternative derivation of Eq. \ref{PC}}

Recalling the implicit information contained in $Z$, we can 
expand the prior for cluster membership $P(Z|\I)$ as, 
\begin{equation}\label{PZI}
\begin{array}{rl}
P(Z|\I) &= P(Z,N_1(Z), ... , N_M(Z), M(Z) | \I )
\\
&=P(Z,N_1(Z), ... , N_m(Z)|M=m, \I)
P(M=m|\I)
\end{array}
\end{equation}
Assuming all assignments $Z=\{Z_i=z_i\}$ are equally likely, then, 
\begin{equation}\label{distP}
\begin{array}{rl}
P(Z,N_1(Z), ... , N_m(Z) | M=m, \I ) 
&= \frac{\#( Z : N_1, ... , N_m)}{\#(Z : N_1...N_m \geq 1) }
\\
&=\frac{n!}{N_1! ... N_m!} \frac{1}{T(n,m)}
\end{array}
\end{equation}
where $T(n,m)$ is the number of distributions of $n$ identifiable items 
into $m$ identifiable boxes \cite{Marcus}. 
Note that this is different to the multinomial distribution with equal 
probabilities $\pi=1/m$ for bin occupancy, because that allows empty 
bins with $N_i=0$, whereas $\{N_i\}$ are counting the number of assignments of 
$Z_i$ to a cluster $g$. 
Eq. \ref{distP} is the probability for a distribution of $N$ elements into $1...m$ bins 
such that no bins are empty, and all independent 
assignments $Z=\{Z_i=z_i\}$ are equally likely. 
Combining Eqs. \ref{PZI}, \ref{distP} and \ref{LZ}, gives, 
\begin{equation}\label{PZ}
\begin{array}{rl}
P(Z|X,\I) 
\propto 
&
\frac{n!}{N_1! ... N_m!}
\frac{1}{T(n,m)}
\\
&\times 
\Pi_{g=1}^m 
\int_{-\infty}^{\infty} d\mu_g f(\mu_g ; \mu_0, \Gamma_0 )  
\left(  \Pi_{i\in C_g}  f\left( x_i ; \mu_g, \Gamma_i \right)  \right)
\end{array}
\end{equation}
Now let $r_j$ be the number of clusters of size $j$, with $j=1..k$, and 
$r_j$ can be zero. 
Then as noted in the main text, there are $n!/N_1! ... N_m!r_1! ... r_k!$  
partitions of type $(N_1, ... , N_m)$. 
A partition is a set of unlabelled clusters, with a form of partial 
labelling implied by the number of elements in the clusters. 
The number of equivalent rearrangements of $m$ unlabelled clusters, with $\{r_j\}$ 
clusters with size $j$, and $m=\sum_{j=1}^k r_j$, is \cite{Marcus}, 
\begin{equation}
\left(
\begin{array}{c}
m=\sum_{j=1}^k r_j
\\
r_1 ... r_k
\end{array}
\right)
\end{equation}
where as before, $r_j$ can be zero.
Each arrangement corresponds to a distribution described by 
Eq. \ref{PZ}, that is unchanged by a permutation of 
the cluster labels.
Therefore, 
\begin{equation}\label{PC2}
\begin{array}{rl}
P(C|\I) & = P(Z|\I)\left(
\begin{array}{c}
m
\\
r_1 ... r_k
\end{array}
\right)
\\
& =  \frac{n!}{N_1! ... N_m!r_1! ... r_k!} \frac{1}{S(n,m)} P(M=m|\I)
\end{array}
\end{equation}
where we used Eq. \ref{distP} and $S(n,m)=T(n,m)/m!$ \cite{Marcus}.
Eq. \ref{PC2} is identical to the combination of Eqs. \ref{PC0} and 
\ref{PCZ} of the main text. 

\section{The prior $P(M=m|\I)$}\label{priorPm}

The main text suggested that the prior $P(M=m|\I)$, might best be chosen 
using the prior information for the particular problem being considered. 
Here we explore the form of $P(M=m|\I)$ that would result from 
randomly partitioning $n$ identical items into $m$ 
parts. 
Taking the number of partitions of $n$ identical items into $m$ parts 
as equivalent to the number of partitions of an integer $n$ into $m$ parts,  
and taking all partitions as equally likely, 
\begin{equation}
P(M=m|\I) = \frac{\# 
	\left(\mbox{Partitions of $n$ items into $m$ non-empty parts} \right)}
{\#\left( \mbox{Partitions of $n$ items} \right) }
\end{equation}
This can be approximated by a formula 
due to Paul Erdos and Joseph Lehner, that gives the distribution 
for the number of partitions of $n$ into 
$m$ elements or less \cite{Erdos},  with as $n \rightarrow \infty$, 
\begin{equation}
F_m(n) \rightarrow \exp \left( 
-\frac{2}{C} \exp \left(
-\frac{C}{2} \frac{m}{\sqrt{n}} + \frac{\log(n)}{2}
\right)
\right) 
\end{equation}
with $C=\pi\sqrt{2/3}$. 
Noting that $\exp(\log(n)/2)=\sqrt{n}$, this may be written as, 
\begin{equation}\label{EqFm}
F_m(n) = \exp\left( - \beta \exp(-k/\beta) \right)
\end{equation}
with $\beta=2\sqrt{n}/C$. 
Eq. \ref{EqFm} can in turn be written as, 
\begin{equation}
F_m(n) = \exp\left( - \exp \left(- \frac{(m-\mu)}{\beta} \right) \right)
\end{equation}
with $\mu=\beta\log(\beta)$. 
This is a Gumbel distribution, with mode $\mu=\beta \log(\beta)$ and 
variance $\sigma^2=\pi^2\beta^2/6=n$. 
However, the distribution is not symmetrical in $n$, and more importantly, 
$\sigma/\mu \sim 1/\log(\sqrt(n)) \rightarrow 0$ as $n\rightarrow \infty$, 
indicating that the distribution becomes increasingly sharply peaked about its 
mode as $n\rightarrow \infty$. 
This would suggest that if all partitions are equally likely, then for 
large datasets, the prior strongly influences the number of clusters. 
This would be a surprising result, and needs further consideration. 
However it would certainly be unsuitable for situations such as a meta-analysis, 
where we expect that $m$ is likely to be $1$.
Therefore for the examples here we will take $P(M=m|\I)=1/n$, and 
leave a more principled choice of prior as a topic for further study. 

\section{Sharply peaked prior}\label{sharpPrior}

The limit of a prior that is sharply peaked around $\mu_0$ can be considered by 
writing $\Gamma_0^{-1}=I \sigma_0^2$, ($\Gamma_0=I/\sigma_0^2$), with $I$ the 
identity matrix, and taking the limit $\sigma_0^2 \rightarrow 0$.  
Firstly expand $\tilde{\Gamma}_g^{-1}=( \sum_{i\in C_g} \Gamma_i )^{-1}$ in 
terms of $\sigma_0^2$,  
\begin{equation}
\begin{array}{l}
\left( \sum_{i\in C_g} \Gamma_i \right)^{-1}
\\
=
\left( \frac{1}{\sigma_0^2} I + \sum_{i\in C_g, i\neq 0} \Gamma_i \right)^{-1}
\\
=
\sigma_0^2 
\left( I + \sigma_0^2 \sum_{i\in C_g, i\neq 0} \Gamma_i \right)^{-1}
\\
=
\sigma_0^2 
\left( 
I - 
\sigma_0^2 \sum_{i\in C_g, i\neq 0} \Gamma_i + 
\sigma_0^4 \sum_{i\in C_g, i\neq 0} \Gamma_i \sum_{j\in C_g, j\neq 0} \Gamma_j
+  \dots \right)
\end{array}
\end{equation}
Then expand $\tilde{\mu}_g$ in terms of $\sigma_0^2$,  
\begin{equation}\label{tildeMuApprox}
\begin{array}{ll}
\tilde{\mu_g} 
&= 
\left( \sum_{i\in C_g} \Gamma_i \right)^{-1}
\sum_{j\in C_g} \Gamma_j  x_j
\\
&= \sigma_0^2 \left( 
I - 
\sigma_0^2 \sum_{i\in C_g, i\neq 0} \Gamma_i + O\left(\sigma_0^4\right) \right) 
\left( \frac{\mu_0}{\sigma_0^2} + \sum_{j\in C_g, j\neq 0} \Gamma_j x_j
\right)
\\
&= \mu_0 +\sigma_0^2  
\sum_{i\in C_g, i\neq 0} \Gamma_i \left( x_i - \mu_0 \right) +  O\left(\sigma_0^4\right) 
\end{array}
\end{equation}
Taking the limit $\sigma_0^2 \rightarrow 0$, 
\begin{equation}\label{limitGamma0}
\begin{array}{l}
\sum_{i \in C_g}
\left( x_i - \tilde{\mu}_g \right)^T \Gamma_i 
\left( x_i - \tilde{\mu}_g \right)
\\
=
\left( \mu_0 - \tilde{\mu}_g \right)^T \Gamma_0 
\left( \mu_0 - \tilde{\mu}_g \right)
+
\sum_{i\in C_g, i\neq 0} 
\left( x_i - \tilde{\mu}_g \right)^T \Gamma_i 
\left( x_i - \tilde{\mu}_g \right)
\\
\rightarrow 
\sum_{i\in C_g, i\neq 0} 
\left( x_i - {\mu}_0\right)^T \Gamma_i 
\left( x_i - {\mu}_0 \right)
+ O\left(\sigma_0^2\right) 
\end{array}
\end{equation}
where the last line used Eq. \ref{tildeMuApprox} to substitute for 
$\tilde{\mu}_g$, and $\Gamma_0=I/\sigma_0^2$. 
Using Eq. \ref{meanI}, the 
other term involving $\Gamma_0$ in Eq. \ref{logLik} may be written as, 
\begin{equation} 
\sum_{g=1}^m p \log n_g 
- \sum_{g=1}^m 
\log  
\frac{\left|\bar{\Gamma}_g^{-1} \right|}{\left| \Gamma_0^{-1} \right|} 
=
\sum_{g=1}^m 
\log  
\frac{\left| \Gamma_0^{-1} \right|}{\left|\tilde{\Gamma}_g^{-1} \right|} 
\end{equation}
where $\bar{\Gamma}_g=\frac{1}{n_g}\sum_{i\in C_g} \Gamma_i$ and 
$\tilde{\Gamma}_g=\sum_{i\in C_g} \Gamma_i$.
Using $|AB|=|A||B|$ and $|B^{-1}|=1/|B|$, this can be written as, 
\begin{equation}\label{logDet}
 \sum_{g=1}^m \log \left| \Gamma_0^{-1}  \sum_{i\in C_g} \Gamma_i \right| 
= 
 \sum_{g=1}^m \log \left| I + \sum_{i\in C_g, i\neq 0} 
\Gamma_0^{-1} \Gamma_i \right|
\end{equation}
Noting that $\Gamma_0^{-1}=I \sigma_0^2$, then as $\sigma_0^2 \rightarrow 0$, 
the determinant on the right-side of Eq. \ref{logDet} tends to $1$, and its 
logarithm tends to zero. 
Therefore as $\sigma^2_0 \rightarrow 0$, the only remaining terms that originally 
involved $\Gamma_0$ are the right-side of Eq. \ref{limitGamma0}. 

\section{Bhattacharyya distance}\label{BD}

The Bhattacharyya distance between two probability densities $p_1(x)$ 
and $p_2(x)$ is, 
\begin{equation}\label{DBCdef}
D_{BC} = 
\int_{-\infty}^{\infty} dx_1 ...
\int_{-\infty}^{\infty} dx_p 
\sqrt{ p_1(x) p_2(x) }
\end{equation}
For two multivariate normals as in Eq. \ref{toFact}, with $\mu_g$ replaced by $x$, 
$x_1$ replaced by $\mu_1$, and $x_2$ replaced by $\mu_2$, this can be integrated 
analytically. 
Using Eq. \ref{ss2} with $C_g=\{1,2\}$ and $\mu$ replaced by $x$, then 
evaluating the normal integral involving $x$ gives, 
\begin{equation}
\begin{array}{c}
\int_{-\infty}^{\infty} dx_1 ...
\int_{-\infty}^{\infty} dx_p 
p_1(x) p_2(x) \\
= \sqrt{ \frac{|2\pi(\Gamma_1+\Gamma_2)^{-1}|}
	{|2\pi\Gamma_1^{-1}||2\pi\Gamma_2^{-1}|} } 
\exp \left\{ 
-\frac{1}{2} \left( \mu_1 - \mu_2 \right)^T 
\Gamma_1 \left(\Gamma_1+\Gamma_2\right)^{-1} \Gamma_2 
\left( \mu_1-\mu_2 \right)
\right\}
\end{array}
\end{equation}
Noting that, 
\begin{equation}
\begin{array}{ll} 
\Gamma_1 \left( \Gamma_1+\Gamma_2 \right)^{-1} \Gamma_2 
&= \left( \Gamma_2^{-1} (\Gamma_1+\Gamma_2) \Gamma_1^{-1} \right)^{-1} 
\\
&= \left( \Gamma_2^{-1} \Gamma_1 \Gamma_1^{-1} 
+ \Gamma_2^{-1} \Gamma_2 \Gamma_1^{-1} \right)^{-1} 
\\
&= \left( \Gamma_1^{-1} + \Gamma_2^{-1} \right)^{-1} 
\end{array}
\end{equation}
and repeating the calculation with the modifications needed to incorporate 
the square roots from the definition (Eq. \ref{DBCdef}), 
\begin{equation}
\begin{array}{c}
\int_{-\infty}^{\infty} dx_1 ...
\int_{-\infty}^{\infty} dx_p 
\sqrt{p_1(x) p_2(x)} 
\\
= \sqrt{\frac{|(\Gamma_1+\Gamma_2)^{-1}|}
	{\sqrt{|\Gamma_1^{-1}||\Gamma_2^{-1}|}} }
\exp \left\{ 
-\frac{1}{8} \left( \mu_1 - \mu_2 \right)^T 
\left(\frac{\Gamma_1^{-1}+\Gamma_2^{-1}}{2}\right)^{-1} 
\left( \mu_1-\mu_2 \right)
\right\}
\end{array}
\end{equation}
where with the square root from Eq. \ref{DBCdef}, the factors of $2\pi$ cancel. 
Note that a factor of $1/2$ has been incorporated into 
$\left(\frac{\Gamma_1^{-1}+\Gamma_2^{-1}}{2}\right)^{-1}$. 
Replacing $\Gamma_i^{-1}$ with $\Sigma_i$ gives the 
Bhattacharyya distance between two multivariate normal distributions in 
its usual form.

\section*{Data and code availability}

Summary data \cite{Webster2021}, test data, and code used to 
produce the figures and tables 
will be made publicly available after publication. 

\section*{Acknowledgements}

AJW thanks Geoff Nicholls for helpful discussions about 
Section \ref{posterior} and advice about 
previous versions of this article. 
Anthony Webster was supported by a fellowship from the 
Nuffield Department of Population Health, University of Oxford, UK.
The funders had no role in study design, data collection and analysis, decision to publish, or preparation of the manuscript.
The results arose from research detailed in Webster et al.\cite{Webster2021}, that 
was conducted using data from UK Biobank,  
a major biomedical database,  under application number 42583. 




\end{document}